\documentclass{article}
\usepackage{amsmath,amsthm,amscd,amssymb,amsfonts}
\input{amssym.def}
\input{amssym.tex}

\newtheorem{thm}{Theorem}[section]
\newtheorem{prop}{Proposition}[section]
\newtheorem{lm}{Lemma}[section]

\newtheorem{fact}{Fact}[section]
\newtheorem{cor}{Corollary}[section]

\title{On a special value of the Ruelle L-function}
\author{Ken-ichi SUGIYAMA
\footnote{Address : Ken-ichi SUGIYAMA,
 Department of Mathematics and Informatics,
Faculty of Science,
Chiba University,
1-33 Yayoi-cho Inage-ku,
Chiba 263-8522, Japan}
\footnote{e-mail address : sugiyama@math.s.chiba-u.ac.jp}}
\begin{document}
\maketitle
\begin{abstract}
Let $X$ be a  a complete hyperbolic threefold of a
 finite volume with only one cusp.
For a unitary local system $\rho$ of rank one on $X$, one may associte
 the Ruelle L-function $R_{\rho}(z)$. Suppose the restriction of $\rho$ to the cusp is
 nontrivial. We will show that the Ruelle L-function has a pole at the
 origin whose order is equal to $-2\dim H^1(X,\,\rho)$. Moreover we will
 prove if
 $\dim H^1(X,\,\rho)$ is zero  $R_{\rho}(0)$ is equal to the square of the Franz-Reidemeister torsion
 of $(X,\,\rho)$.
\end{abstract}
\section{Introduction}
In \cite{Sugiyama2} we have shown that a geometric analog of the Iwasawa
conjecture holds for the Ruelle L-function and the twisted Alexander
invariant.\\

More precisely let $\Gamma$ be a torsion free cofinite discrete subgroup of
$PSL_{2}({\mathbb C})$. It acts on the three dimensional Poincar\'{e}
upper half space
\[
 {\mathbb H}^3=\{(x,\,y,\,r)\,|\,x,\,y\in{\mathbb R},\,r>0\}
\]
endowed with a metric
\[
 ds^2=\frac{dx^2+dy^2+dr^2}{r^2},
\]
whose sectional curvature $\equiv -1$.
 Let $X$ be the quotient, which is a complete hyperbolic threefold of finite volume. We will
assume that it has only one cusp. 
Let $\rho$ be a unitary character of $\Gamma$. It defines a unitary
local system on $X$ of rank one, which will be denoted by the same symbol.
 By the one to one
correspondence between the set of loxiodromic conjugacy classes of
$\Gamma$ and one of closed geodesics of $X$,  {\it the Ruelle
L-function} is defined as
\[
 R_{\rho}(z)=\prod_{\gamma}\det[1-\rho(\gamma)e^{-zl(\gamma)}],
\]
where $\gamma$ runs through primitive closed
geodesics. Here $z$ is a complex number and $l(\gamma)$ is the length of
$\gamma$. 
It is known $R_{\rho}(z)$ is absolutely convergent if ${\rm Re}\,z$ is sufficiently large. Suppose the
restriction $\rho|_{\Gamma_{\infty}}$ of $\rho$ to the fundamental group
$\Gamma_{\infty}$ of the cusp is nontrivial. In \cite{Sugiyama2} we have shown that $R_{\rho}(z)$
is meromorphically continued on the whole plane and that
\[
 {\rm ord}_{z=0}R_{\rho}(z)=-2h^{1}(\rho),
\]
where $h^{1}(\rho)$ is the dimension of $H^{1}(X,\,\rho)$.\\

Let us assume there is a surjective homomorphism from $\Gamma$ to ${\mathbb
Z}$ and $X_{\infty}$ the corresponding infinite cyclic covering of
$X$. Moreover suppose that the dimensions of all of
$H_{\cdot}(X_{\infty},\,{\mathbb C})$ and $H_{\cdot}(X_{\infty},\,\rho)$
are finite. Let $g$ be a generator of the infinite cyclic group. Then the twisted Alexander invariant $A_{X}^{*}(\rho)$ is
defined to be an alternating product of characteristic polynomials of
the action of $g$ on $H^{\cdot}(X_{\infty},\,\rho)$. (See
\cite{Sugiyama} for the precise definition.) In \cite{Sugiyama} we have prove that
\[
 {\rm ord}_{z=0}R_{\rho}(z)\geq 2{\rm ord}_{t=1}A_{X}^{*}(\rho),
\]
and that if $h^{1}(\rho)$ is zero, $R_{\rho}(z)$ and
$A_{X}^{*}(\rho)(t)$ does not vanish at $z=0$ and $t=1$,
respectively. It should be natural to compare their values. In
fact we will prove the following theorem. 

\begin{thm} Suppose that $\rho|_{\Gamma_{\infty}}$ is nontrivial and
 that $h^1(\rho)$ vanishes. Then we have
\[
 R_{\rho}(0)=\tau_X(\rho)^2,
\]
where $\tau_X(\rho)$ is the Reidemeister torsion of $X$ and $\rho$.
\end{thm}
If the manifold is compact, the corresponding result has been already proved by
Fried (\cite{Fried}).\\

Since we know the absolute value of $\tau_X(\rho)$ is equal to
a product of $|A_{X}^{*}(\rho)(1)|$ and a certain positive constant $\delta_{\rho}$ which
can be computed explicitly (\cite{Sugiyama} {\bf Theorem 3.4}), we have
\[
 |R_{\rho}(0)|=(\delta_{\rho}\cdot |A_{X}^{*}(\rho)(1)|)^2.
\]

{\bf Acknowledgement.}\hspace{5mm}
It is a great pleasure to appreciate Professor Park for his kindness to
answer our many questions, as well as Professor Wakayama for sending his
manuscripts which were great hepl for us. It is clear without their help
our work will not be completed.

\section{Laplace-Mellin transform}
We define the Laplace transform of a function $f$ on ${\mathbb R}$ to be
\[
 L(f)(z)=\int_{0}^{\infty}e^{-tz^2}\frac{f(t)}{t}dt,
\]
if the RHS is absolutely convergent.
\begin{lm}
Let $l$ be a positive number and suppose $z>0$. Then
\[
 L(\frac{1}{\sqrt{4\pi t}}e^{-\frac{l^2}{4t}})(z)=\frac{e^{-lz}}{l}.
\]
\end{lm}
{\bf Proof.}
Taking a derivative of 
\[
 \int^{\infty}_0\exp(-t^2-\frac{x^2}{t^2})dt=\frac{\sqrt{\pi}}{2}e^{-2x}
\]
with respect to $x$, we have
\[
 x\int^{\infty}_0\frac{1}{t^2}\exp(-t^2-\frac{x^2}{t^2})dt=\frac{\sqrt{\pi}}{2}e^{-2x}.
\]
Let $\alpha$ be a positive number. If we make a change of variables:
\[
 t \to \sqrt{\alpha t},
\]
we will obtain
\begin{equation}
\int^{\infty}_0t^{-\frac{3}{2}}\exp(-t^2-\frac{x^2}{t^2})dt=\frac{\sqrt{\pi \alpha}e^{-2x}}{x}.
\end{equation}
Now (1) and a simple computation will show the desired identity.
\begin{flushright}
$\Box$
\end{flushright}
We also define the Laplace-Mellin transform of $f$ to be
\[
 {\mathcal L}(f)(s,\,z)=\int^{\infty}_{0}e^{-tz^2}t^{s-1}f(t)dt, 
\]
for sufficiently large real numbers $z$ and $s$ if the RHS
is absolutely convergent. Suppose that ${\mathcal L}(f)(s,\,z)$ is continued
to a meromorphic function on an open domain $U$ of ${\mathbb C}^2$ which
contains
\[
 \{(s,\,z)\,|\, s,\,z\in{\mathbb R}, \,s,\,z>>0,\},
\]
 and that its polar set $P_{{\mathcal L}(f)(s,\,z)}$ does not contain
\[
 U_{0,z}=U\cap {\mathbb C}_{0,z},
\]
where ${\mathbb C}_{0,z}=\{(0,\,z)\,|\, z\in {\mathbb C}\}$. Then we define
the Laplace transform $L(f)(z)$ on $U_{0,z}$ to be
\[
 L(f)(z)={\mathcal L}(f)(0,\,z).
\]
\\

For a nonnegative integer $k$, let us consider a function:
\[
 p_{k}(t)=\int^{\infty}_0e^{-tx^2}x^{2k}dx.
\] 
\begin{lm}
For $z>0$ and $s>\frac{1}{2}+k$, the Laplace-Mellin transform of $p_k$ is
\[
 {\mathcal L}(p_k)(s,\,z)=\frac{\sqrt{\pi}C_k}{2}z^{1+2k-2s}\Gamma(s-\frac{1}{2}-k),
\]
which is defined over $\{(s,\,z)\in{\mathbb C}^2\,|\, s\in{\mathbb C},\,-\pi<{\rm Im}\,z< \pi\}$.
Here we put
\[
 C_{0}=1
\]
and
\[
 C_{k}=\prod_{m=0}^{k-1}(m+\frac{1}{2})
\]
for $k\geq 1$.
\end{lm}
{\bf Proof.}
Let $t$ be a positive number. Take the $k$-times derivative of
\[
 \int^{\infty}_{-\infty}e^{-tx^2}dx=\sqrt{\pi}t^{-\frac{1}{2}}
\]
with respect to $t$, we obtain
\begin{equation}
 \int^{\infty}_{-\infty}x^{2k}e^{-tx^2}dx=\frac{\sqrt{\pi}C_k}{2}t^{-\frac{1}{2}-k}.
\end{equation}
Then we compute:
\begin{eqnarray*}
 {\mathcal L}(p_k)(s,\,z)&=&\int^{\infty}_{0}\frac{dt}{t}t^{s}e^{-tz^2}\int^{\infty}_{-\infty}x^{2k}e^{-tx^2}dx\\
&=&\frac{\sqrt{\pi}C_k}{2}\int^{\infty}_{0}e^{-tz^2}t^{s-\frac{1}{2}-k}\frac{dt}{t}\\
&=&\frac{\sqrt{\pi}C_k}{2}z^{1+2k-2s}\Gamma(s-\frac{1}{2}-k).
\end{eqnarray*}
\begin{flushright}
$\Box$
\end{flushright}

\begin{cor}
For a nonnegative integer $k$,
\[
 L(p_k)(z)=\frac{\sqrt{\pi}C_k}{2}\Gamma(-\frac{1}{2}-k)z^{1+2k}.
\]
\end{cor}
Note that this is defined over the whole $z$-plane.
\section{Selberg trace formula and Laplace transforms of orbital
 integrals}

Let $\Omega_{X}^j(\rho)$ be a vector bundle of $j$-forms on $X$
twisted by $\rho$ and the space of its square integrable sections will be
denoted by $L^2(X,\,\Omega_{X}^j(\rho))$. The positive Hodge Laplacian has the selfadjoint extension to $L^2(X,\,\Omega_{X}^j(\rho))$,
which will be denoted by $\Delta$. Note that the Hodge star operator
induces an isomorphism of Hilbert spaces:
\begin{equation}
 L^2(X,\,\Omega_{X}^j(\rho))\simeq L^2(X,\,\Omega_{X}^{3-j}(\rho)),\quad j=0,\,1,
\end{equation}
which commutes with $\Delta$. \\
Since $\rho|_{\Gamma_{\infty}}$ is nontrivial we know that the spectrum of
$\Delta$ consists of only eigenvalues and the Selberg trace formula
becomes (See $\S4$ of \cite{Sugiyama2}):
\[
 {\rm Tr}[e^{-t\Delta}\,|\,L^2(X,\,\Omega_{X}^j(\rho))]={\mathcal
 I}_{j}(t)+{\mathcal H}_{j}(t)+{\mathcal U}_{j}(t).
\]
Here ${\mathcal I}_{j}(t)$,  ${\mathcal H}_{j}(t)$ and  ${\mathcal
U}_{j}(t)$ are the identity, the hyperbolic and the unipotent term,
respectively. In this section we will compute their Laplace transform.
\begin{enumerate}
\item {\bf The hyperbolic term}
 
Let $A$ be a split Cartan subgoup of $G=PSL_{2}({\mathbb C})$. The Lie
algebras of $G$ and $A$ will be denoted by ${\frak G}$ and ${\frak A}$, respectively. The choice of $A$ determines a positive root $\alpha$ of ${\frak G}$ and let $H$ be an
element of ${\frak A}$ satisfying
\[
 \alpha(H)=1.
\]
If we exponentiate a linear isomorphism:
\[
 {\mathbb R} \stackrel{h}\to {\frak A},\quad h(t)=tH, 
\]
we know $A$ is isomorphic to the multiplicative group of positive real
numbers ${\mathbb R}_{+}$ and will identify them.

Let $K\simeq SO_{3}({\mathbb R})$ be the maximal compact
subgroup. According to the Iwasawa decompostion $G=KAN$ an element $g$
of $G$ can be written as
\[
 g=k(g)a(g)n(g).
\]

Let $M$ be the centralizer of $A$ in
$K$, which is isomorphic to $SO_{2}({\mathbb R})$. It determines a
paraboloic subgroup: 
\begin{equation}
 P=MAN.
\end{equation}
Let $\Gamma_h$ be the set of
conjugacy classes of loxidromic elements of $\Gamma$. Since there is a
natural bijection between closed geodesics of $X$ and $\Gamma_h$, we may
identify them. Thus an element $\gamma$ of $\Gamma_h$ is written as 
\[
 \gamma=\gamma_0^{\mu(\gamma)},
\]
where $\gamma_0$ is a primitive closed geodesic and $\mu(\gamma)$ is a
positive integer, which will be referred as {\it the multiplicity}. 
The length of $\gamma\in \Gamma_h$ will be denoted by $l(\gamma)$ and 
let $\Gamma_{h,prim}$ be the set of primitive closed geodesics.\\

Using the Langlands decomposition (4), $\gamma\in \Gamma_h$ may be written as
\[
 g\gamma g^{-1}=m(\gamma)\cdot a(\gamma)\in MA
\]
for a certain $g\in G$. Here $m(\gamma)$ is nothing but the
holonomy of a pararell transformation along $\gamma$. Note that elements of $GL_2({\mathbb R})$:
\[
 A^{u}(\gamma)=e^{l(\gamma)}m(\gamma), \quad A^{s}(\gamma)=e^{-l(\gamma)}m(\gamma)
\]
describe an unstable or a stable action of the linear Poincar\'{e}
map, respectively.\\

 For $\gamma\in \Gamma_h$ we set
\[
 \Delta(\gamma)=\det [I_2-A^{s}(\gamma)]
\]
and
\[
 a_{0}(\gamma)=\frac{\rho(\gamma)\cdot l(\gamma_0)}{\Delta(\gamma)},\quad
 a_{1}(\gamma)=\frac{\rho(\gamma)\cdot {\rm Tr}\,[m(\gamma)]\cdot l(\gamma_0)}{\Delta(\gamma)}.
\]
Now {\bf Theorem 2} of \cite{Fried} shows the hyperbolic terms are
given by
\[
 {\mathcal H}_0(t)=H_{0}(t),\quad {\mathcal H}_1(t)=H_{0}(t)+H_{1}(t),
\]
where 
\[
 H_{0}(t)=\Sigma_{\gamma\in \Gamma_h}\frac{a_0(\gamma)}{\sqrt{4\pi
 t}}\exp [-(\frac{l(\gamma)^2}{4t}+t+l(\gamma))],
\]
and
\[
 H_{1}(t)=\Sigma_{\gamma\in \Gamma_h}\frac{a_1(\gamma)}{\sqrt{4\pi
 t}}\exp [-(\frac{l(\gamma)^2}{4t}+l(\gamma))].
\]
We will explain a relation between these hyperbolic terms and the
Ruelle L-function\\

 For $j=0,\,1$ we set
\[
 S_j(z)=\exp [-\sum_{\gamma\in \Gamma_{h}}\frac{a_j(\gamma)}{l(\gamma)}e^{-zl(\gamma)}].
\]
Then the formula (RS) of \cite{Fried} shows
\[
 R_{\rho}(z)=\frac{S_{0}(z)S_{0}(z+2)}{S_{1}(z+1)}.
\]
Using {\bf Lemma 2.1}, a simple computation implies the
      following lemma. 
\begin{lm}
\begin{enumerate}
\item
\[
 L(H_1)(z)=-\log S_{1}(z+1).
\]
\item
\[
 L(e^tH_0)(z)=-\log S_{0}(z+1).
\]
\end{enumerate}
\end{lm}

Thus we have proved the following proposition.
\begin{prop}
\[
 \log R_{\rho}(0)=L(H_1)(0)-L(e^tH_0)(-1)-L(e^tH_0)(1)
\]
\end{prop}

\item {\bf The identity term}
In $\S6$ of \cite{Sugiyama2} we have computed the identity terms to be:
\[
 {\mathcal I}_0(t)=I_{0}(t),\quad {\mathcal I}_1(t)=I_{0}(t)+I_{1}(t),
\]
where 
\[
 I_{0}(t)=vol(X)\int^{\infty}_{-\infty}e^{-t(x^2+1)}x^2dx,
\]
and
\[
 I_{1}(t)=2vol(X)\int^{\infty}_{-\infty}e^{-tx^2}(x^2+1)dx.
\]

{\bf Lemma 2.2} implies
\[
 {\mathcal L}(e^tI_0)(s,\,z)=\frac{\sqrt{\pi}}{4}vol(X)z^{3-2s}\Gamma(s-\frac{3}{2}).
\]
Also the identity
\[
 \Gamma(-\frac{3}{2})=\frac{4\sqrt{\pi}}{3}
\]
implies 
\[
 L(e^tI_0)(z)=\frac{\pi}{3}vol(X)z^{3}.
\]
By the same computation, we will have
\[
 {\mathcal L}(I_1)(s,\,z)=\frac{\sqrt{\pi}}{2}vol(X)(z^{3-2s}\Gamma(s-\frac{3}{2})+2z^{1-2s}\Gamma(s-\frac{1}{2})),
\]
and
\[
 L(I_1)(z)=2\pi vol(X)(\frac{z^{3}}{3}-z).
\]
Thus we have proved
\begin{prop}
\[
 L(I_1)(0)-L(e^tI_0)(-1)-L(e^tI_0)(1)=0.
\]
\end{prop}

\item {\bf The unipotent term}\\

We put
\[
 U_0(t)={\mathcal U}_0(t),\quad U_1(t)={\mathcal U}_1(t)-{\mathcal U}_0(t).
\]
In {\bf Proposition 7.1} of \cite{Sugiyama2}, we have proved the following fact.
\begin{fact}
\begin{enumerate}
\item 
\[
 U_0(t)=C_{\rho,\Gamma}e^{-t}\int_{-\infty}^{\infty}e^{-tx^2}dx,
\]
\item
\[
 U_{1}(t)=2C_{\rho,\Gamma}\int_{-\infty}^{\infty}e^{-tx^2}dx.
\]
\end{enumerate}
where $C_{\rho,\Gamma}$ is a constant determined by $\Gamma$ and $\rho$.
\end{fact}
Thus by {\bf Lemma 2.2}, we obtain
\begin{eqnarray*}
{\mathcal L}(U_{1})(s,\,z)&=& 2{\mathcal L}(e^tU_{0})(s,\,z)\\
&=& \sqrt{\pi}C_{\rho,\Gamma}z^{1-2s}\Gamma(s-\frac{1}{2}),
\end{eqnarray*}
which implies
\begin{eqnarray*}
L(U_{1})(z)&=& 2L(e^tU_{0})(z)\\
&=& -\pi C_{\rho,\Gamma}z.
\end{eqnarray*}
Thus the following proposition is proved.
\begin{prop}
\[
 L(U_{1})(0)-L(e^tU_{0})(-1)-L(e^tU_{0})(1)=0.
\]
\end{prop}
\end{enumerate}
\section{Laplace transform of the heat kernel and the analytic torsion}
We set
\[
 \delta_{0,\rho}(t)={\rm Trace}[e^{-t\Delta}\,|\,L^{2}(X,\,\Omega_{X}^{0}(\rho))],
\]
and 
\[
  \delta_{1,\rho}(t)={\rm Trace}[e^{-t\Delta}\,|\,L^{2}(X,\,\Omega_{X}^{1}(\rho))]-\delta_{0,\rho}(t).
\]
The nontriviality of $\rho|_{\Gamma_{\infty}}$ implies $H^0(X,\,\rho)=0$ and by
the Zucker's result (\cite{Zucker}, see also the
introduction of \cite{Mazzeo-Phillips} and $\S2$ of \cite{Sugiyama2}), we have
\[
 {\rm Ker}\,[\Delta\,|\,L^{2}(X,\,\Omega_{X}^{0}(\rho))]=0.
\]
Let us assume $h^1(\rho)$ vanishes. As we have shown \cite{Sugiyama2}
{\bf Lemma 2.1}, this implies
\[
 {\rm Ker}\,[\Delta\,|\,L^{2}(X,\,\Omega_{X}^{1}(\rho))]=0.
\]
Thus there is positive constants $c_j$ and $A$ such that
\begin{equation}
 |\delta_j(t)|\leq c_0e^{-c_1t^2}\quad \mbox{for}\quad t\geq A.
\end{equation}
\begin{lm}
\begin{enumerate}
\item ${\mathcal L}(\delta_1)(s,\,z)$ is absolutely
      convergent for ${\rm Re}\,s>>0$ and $z>0$ and is meromorphically
      continued on an open domain of ${\mathbb C}^2$ containing 
\[
 \{(s,\,z)\,|\, s\in {\mathbb C},\,z\in{\mathbb R}\}.
\]
\item 
${\mathcal L}(\delta_0)(s,\,z)$ is absolutely
      convergent for ${\rm Re}\,s>>0$ and $z\geq 1$ is meromorphically
      continued on an open subset of ${\mathbb C}^2$ containing 
\[
 \{(s,\,z)\,|\, s\in {\mathbb C},\,z\geq 1\}.
\]
\end{enumerate}
\end{lm}
{\bf Proof.} Since a proof of the both statements are same, we will only
prove the first. The absolutely convergence is clear from (4).\\

Let us write
\[
 {\mathcal L}(\delta_1)(s,\,z)={\mathcal L}_{(0,A]}(\delta_1)(s,\,z)+{\mathcal L}_{[A,\infty)}(\delta_1)(s,\,z),
\]
where we put
\[
 {\mathcal L}_{(0,A]}(\delta_1)(s,\,z)=\int^{A}_{0}e^{-tz^2}t^{s-1}\delta_1(t)dt,
\] 
and
\[
 {\mathcal L}_{[A,\infty)}(\delta_1)(s,\,z)=\int^{\infty}_{A}e^{-tz^2}t^{s-1}\delta_1(t)dt. 
\] 
(5) implies ${\mathcal L}_{[A,\infty)}(\delta_1)(s,\,z)$ is defined on
    such an open subset. The computation of the previous section and the
    equation (2) show
    the orbital integrals have the following asymptotic expansion when $t\to 0$:
\begin{eqnarray*}
 H_{1}(t)&=&\Sigma_{\gamma\in \Gamma_h}\frac{a_1(\gamma)}{\sqrt{4\pi
 t}}\exp [-(\frac{l(\gamma)^2}{4t}+l(\gamma))]\\
&\sim& \alpha_0e^{-\frac{\alpha_1}{t}},
\end{eqnarray*}
\begin{eqnarray*}
I_{1}(t)&=&2vol(X)\int^{\infty}_{-\infty}e^{-tx^2}(x^2+1)dx\\
&\sim& \beta_1t^{-\frac{3}{2}}+\beta_0t^{-\frac{1}{2}},
\end{eqnarray*}
and
\begin{eqnarray*}
 U_{1}(t)&=&2C_{\rho,\Gamma}\int_{-\infty}^{\infty}e^{-tx^2}dx\\
&\sim& \gamma_0t^{-\frac{1}{2}}.
\end{eqnarray*}
Thus for a real number $z$, using the Selberg trace formula, we have an asymptotic expansion:
\begin{eqnarray*}
{\mathcal L}_{(0,A]}(\delta_1)(s,\,z)&\sim&\alpha_0\int^{A}_0e^{-\frac{\alpha_1}{t}}e^{-tz^2}t^{s-1}dt\\
&+&
 \beta_1\int^{A}_0e^{-tz^2}t^{s-\frac{5}{2}}dt\\
&+& \gamma^{\prime}_0\int^{A}_0e^{-tz^2}t^{s-\frac{3}{2}}dt\\
&\sim& A_0 + \frac{A_1}{s-\frac{5}{2}}+\frac{A_2}{s-\frac{3}{2}}.
\end{eqnarray*}
where $\alpha_0$, $\beta_i$, $\beta^{\prime}_i$, $\gamma_i$,
$\gamma^{\prime}_i$ and $A_i$ are constants.  Now we have obtained the
desired result.
\begin{flushright}
$\Box$
\end{flushright}
If ${\rm Re}\,s$ is sufficiently large, the integral
\[
 {\mathcal L}(\delta_1)(s,\,0)=\int_{0}^{\infty}\delta_1(t)t^{s-1}dt
\]
is absolutely convergent and is nothing but the Mellin transform
$M(\delta_1)(s)$ of $\delta_1$. Since by {\bf Lemme 3.1} of
\cite{Sugiyama2} we know $L(\delta_1)$ is regular at $z=0$, {\bf Lemma
4.1} implies
\[
 L(\delta_1)(0)={\mathcal L}(\delta_1)(0,\,0)=M(\delta_1)(0).
\]
Using {\bf Lemme 3.2} of \cite{Sugiyama2}, the same argument will imply
\[
 L(e^t\delta_0)(1)=M(\delta_0)(0).
\]
In order to compute $L(e^t\delta_0)(-1)$, we prepare the following lemma.
\begin{lm} Let us put
\[
 L_0(z)=L(e^t\delta_1)(z-1).
\]
Then it satisfies a functional equation:
\[
 L_0(1+z)=L_0(1-z).
\]
\end{lm}
{\bf Proof.}
Let $F$ be their difference:
\[
 F(z)=L_0(1+z)-L_0(1-z).
\]
{\bf Lemma 3.2} of \cite{Sugiyama2} shows
\[
 F^{\prime}(z)=L^{\prime}_0(1+z)+L^{\prime}_0(1-z)=0,
\]
and therefore $F$ is a constant. But since 
\[
 \lim_{z\to +\infty}L_0(z)=\lim_{z\to -\infty}L_0(z)=0
\]
we know $F=0$.
\begin{flushright}
$\Box$
\end{flushright}
In particular we have
\[
 L_0(2)=L_0(0),
\]
which implies
\[
 L(e^t\delta_0)(1)=L(e^t\delta_0)(-1).
\]

Thus we have proved the equation: 
\begin{equation}
M(\delta_1)(0)-2M(\delta_0)(0)=L(\delta_1)(0)-L(e^t\delta_0)(-1)-L(e^t\delta_0)(1).
\end{equation}
Using {\bf Proposition 3.1}, {\bf Proposition 3.2} and {\bf Proposition 3.3},
the Selberg trace formula informs us the RHS is equal to $\log
R_{\rho}(0)$.
Thus we have obtained 
\begin{equation}
\log R_{\rho}(0)=M(\delta_1)(0)-2M(\delta_0)(0).
\end{equation}
Now let us recall the definition of the analytic torsion
$T_{X}(\rho)$ of $(X,\,\rho)$ due to Ray and Singer
\cite{Ray-Singer} (See also \cite{Cheeger} and \cite{Muller}):
\[
 \log
 T_{X}(\rho)=\frac{1}{2}\frac{d}{ds}[\frac{1}{\Gamma(s)}\sum_{j=0}^{3}(-1)^{j+1}j\cdot M({\rm Trace}[e^{-t\Delta_{X}}\,|\,L^{2}(X,\,\Omega^{j}(\rho))])(s)]|_{s=0}
\]
As we have seen the Mellin transform of the heat kernel on
$L^{2}(X,\,\Omega^{j}(\rho)),\,(j=0,1)$ is regular at the origin and (3) implies 

\[
 \log T_{X}(\rho)=\frac{1}{2}(2M(\delta_0)(0)-M(\delta_1)(0)).
\]
Thus we have obtained the following theorem.
\begin{thm} Suppose $h^1(\rho)$ vanishes. Then
\[
 R_{\rho}(0)=T_{X}(\rho)^2.
\]
\end{thm}
\section{The theorem of Cheeger and M\"{u}ller}
For a positive number $A$ let $X_A$ the image of
\[
 {\mathbb H}^3_{A}=\{(x,\,y,\,r)\in {\mathbb H}^3\,|\,r\leq A \},
\]
under the natural projection
\[
 {\mathbb H}^3\stackrel{\pi}\to X.
\]
Let $Y_A$ be the complement of the interior of $X_A$.
We take $A$ sufficiently large so that the boundary $M_A$ of $X_A$ is
a flat torus and that $Y_A$ is diffeomorphic to a product of  $M_A$ with an
interval $[A,\,\infty)$. \\


We will review the analytic torsion of $(X_A,\rho)$ with respect to the
absolute boundary condition. Let $\Omega_{X}^{\cdot}(\rho)|_{M_A}$ be th
restricton $\Omega_{X}^{\cdot}(\rho)$ to $M_A$. Its section $\omega$ can
be written as
\[
 \omega=\omega_t+dr\wedge \omega_n,
\]
where $\omega_t$ and $\omega_n$ are sections of
$\Omega_{M_A}^{\cdot}(\rho)$. We put
\[
 P_a(\omega)=\omega_n,
\]
and
\[
 P_r(\omega)=\omega_t.
\]
The space of smooth $j$-forms on $X_A$ (resp. $Y_A$) twisted by $\rho$ satisfying the absolute
(resp. relative) boundary condition is defined to be
\[
 C^{\infty}_{abs}(X_A,\,\Omega^j(\rho))=\{\omega\in
 C^{\infty}(X_A,\,\Omega^j(\rho))\,|\, P_a(\omega)=P_a(d\omega)=0\},
\]
(resp.
\[
 C^{\infty}_{rel}(Y_A,\,\Omega^j(\rho))=\{\omega\in
 C^{\infty}(Y_A,\,\Omega^j(\rho))\,|\, P_r(\omega)=P_r(\delta\omega)=0\},
\]
where $\delta$ is the formal adjoint of $d$.) It is known that elements
of each space satisfy the self-adjoint boundary condition
(\cite{Dai-Fang} (2.8)):
\begin{equation}
\omega,\,\omega^{\prime}\in C^{\infty}_{abs}(X_A,\,\Omega^j(\rho))\Rightarrow (\Delta\omega,\,\omega^{\prime})=(\omega,\,\Delta\omega^{\prime}),
\end{equation}
\begin{equation}
\eta,\,\eta^{\prime}\in C^{\infty}_{rel}(Y_A,\,\Omega^j(\rho))\Rightarrow (\Delta\eta,\,\eta^{\prime})=(\eta,\,\Delta\eta^{\prime}).
\end{equation}
For $\omega\in C^{\infty}_{abs}(X_A,\,\Omega^j(\rho))$ we define
$\tilde{\omega}\in L^2(X,\,\Omega^{j}(\rho))$ to be
\[
\tilde{\omega}(x)=
 \left\{\begin{array}{ccc}
\omega(x)&\mbox{if}&x\in X_A\\
0 &\mbox{if}&x\notin X_A.
\end{array}
\right.
\]
In this way we may consider $C^{\infty}_{abs}(X_A,\,\Omega^j(\rho))$ as
a subspace of $L^2(X,\,\Omega^{j}(\rho))$ and let
$L^{2}_{abs}(X_A,\,\Omega^j(\rho))$ be its closure. By the same
procedure, we have a closed subspace
$L^{2}_{rel}(Y_A,\,\Omega^j(\rho))$.
(8) and (9) implies the positive Laplacian has a selfadjoint extension
$\Delta_{X_A}$ and $\Delta_{Y_A}$ on $L^{2}_{abs}(X_A,\,\Omega^j(\rho))$
and $L^{2}_{rel}(Y_A,\,\Omega^j(\rho))$, respectively. Moreover there is
an orthogonal decomposition:
\[
 L^2(X,\,\Omega^{j}(\rho))=L^{2}_{abs}(X_A,\,\Omega^j(\rho))\hat{\oplus}L^{2}_{rel}(Y_A,\,\Omega^j(\rho)),
\]
which makes $\Delta$ into a block form
\begin{equation}
 \Delta=
\left(
\begin{array}{cc}
\Delta_{X_A,j}&0\\
0&\Delta_{Y_A,j}
\end{array}
\right).
\end{equation}
For a positive $t$, the heat operator $e^{-t\Delta_{X_A,j}}$ is in the trace class
and the integral
\[
 M({\rm Trace}[e^{-t\Delta_{X_A,j}}])(s)=\int^{\infty}_{0}t^{s-1}{\rm Trace}[e^{-t\Delta_{X_A,j}}]dt
\]
is absolutely convergent for ${\rm Re}\,s>>0$. Moreover it is
meromorphically continued on the whole plane. \\

Let us investigate its behavior at the origin. As we have
seen in \cite{Sugiyama2} $\S4$, the
nontriviality of $\rho|_{\Gamma_{\infty}}$ implies the cuspidality of any
element of $L^{2}_{rel}(Y_A,\,\Omega^j(\rho))$. Then the the proof of
\cite{SLNMuller} {\bf Proposition 4.9}  (especially the equation (4.12))
shows the infimum of the set of spectrum of $\Delta_{Y_A,j}$ has the
following lower bound:
\[
 {\rm Inf}\,\sigma(\Delta_{Y_A,j})>cA,
\]
where $c$ is a positive constant. Thus the dimension of the kernel of
$\Delta_{X}$ on $L^2(X,\,\Omega^j(\rho))$ and $\Delta_{X_A,j}$ are same. By the assumption
the previous has the trivial kernel, so does $\Delta_{X_A,j}$. This
implies $M({\rm Trace}[e^{-t\Delta_{X_A,j}}])(s)$ is regular at $s=0$. Now
the analytic torsion $T_{X_A}(\rho)$ of $(X_A,\rho)$ (with respect to to the absolute
boundary condition) is defined to be
\[
 \log
 T_{X_A}(\rho)=\frac{1}{2}\frac{d}{ds}[\frac{1}{\Gamma(s)}\sum_{j=0}^{3}(-1)^{j+1}j\cdot M({\rm Trace}[e^{-t\Delta_{X_A,j}}])(s)]|_{s=0}
\]
Let $\tau_{X_A}(\rho)$ be the Reidemeister torsion of $(X_A,\,\rho)$. If we apply {\bf Theorem 1.1} of \cite{Dai-Fang}, we obtain 
\begin{equation}
 T_{X_A}(\rho)=\tau_{X_A}(\rho).
\end{equation}
Here we have used the following fact. First of all, one can directly
check that the second fundamental form of $M_A$ is zero and therefore
the term $\phi$ in the theorem vanishes. Next since the dimension of $X_A$ is three, the Chern-Simon class
defined by Bisumut and Zhang (\cite{Bismut-Zhang}) also
vanishes. Finally the index theorem inform us the
Euler characteristic $\chi(M_{A},\,\rho)$ is zero. \\

For sufficiently large $A$ and $A^{\prime}$, $X_A$ and
$X_{A^{\prime}}$ are isomorphic as PL-manifolds and the PL-invariance of
the Reidemeister torsion implies
\[
 \tau_{X_A}(\rho)=\tau_{X_{A^{\prime}}}(\rho).
\]
Thus the following definition makes sense:
\begin{equation}
 \tau_{X}(\rho)=\lim_{A\to\infty}\tau_{X_A}(\rho).
\end{equation}
\\

Let us compare the analytic torsion of $(X_A,\,\rho)$ and $(X,\,\rho)$
\begin{prop}
\[
 T_X(\rho)=\lim_{A\to\infty}T_{X_A}(\rho).
\]
\end{prop}
{\bf Proof.}
For ${\rm Re}\,s>>0$, M\"{u}ller's result cited before implies 
\begin{eqnarray*}
M({\rm Trace}[e^{-t\Delta}\,|\,L^{2}(X,\,\Omega^j(\rho))])(s)&=&\int^{\infty}_{0}t^{s-1}{\rm Trace}[e^{-t\Delta}\,|\,L^{2}(X,\,\Omega^j(\rho))]dt\\
&=& \lim_{A\to\infty}\int^{\infty}_{0}t^{s-1}{\rm Trace}[e^{-t\Delta_{X_A,j}}]dt\\
&=& \lim_{A\to\infty}M({\rm Trace}[e^{-t\Delta_{X_A,j}}])(s),
\end{eqnarray*}
which yields the identity as meromorphic functions on the whole plane:
\[
 M({\rm Trace}[e^{-t\Delta}\,|\,L^{2}(X,\,\Omega^j(\rho))])(s)=\lim_{A\to\infty}M({\rm Trace}[e^{-t\Delta_{X_A,j}}])(s).
\]
Now the desired identity will follow from the definition of the analytic torsion.
\begin{flushright}
$\Box$
\end{flushright}
Now {\bf Theorem 4.1}, (11), (12) and {\bf Proposion 5.1} implies the
following theorem.

\begin{thm} Suppose that $\rho|_{\Gamma_{\infty}}$ is nontrivial and
 that $h^1(\rho)$ vanishes. Then
\[
 R_{\rho}(0)=\tau_X(\rho)^2.
\]
\end{thm}


\vspace{10mm}
\begin{flushright}
Address : Department of Mathematics and Informatics\\
Faculty of Science\\
Chiba University\\
1-33 Yayoi-cho Inage-ku\\
Chiba 263-8522, Japan \\
e-mail address : sugiyama@math.s.chiba-u.ac.jp
\end{flushright}
\end{document}